\documentclass[12pt,reqno,a4paper]{article}
\usepackage{filecontents}

	\usepackage[outer=2.5cm,top=2.5cm,bottom=2.5cm,inner=2.0cm]{geometry}

	\usepackage{xfrac}
	\usepackage{amsmath}
	\usepackage[english]{varioref}
	\usepackage[bookmarksopen=true, bookmarksnumbered=true, hidelinks]{hyperref}
	\usepackage[english]{cleveref}\usepackage[T1]{fontenc}
	\usepackage{amsthm}
	\usepackage{amssymb}
	\usepackage{amsfonts}
	\usepackage{dsfont}
	\usepackage{stmaryrd}
	\usepackage{commath}
	\usepackage{mathrsfs}
	\usepackage{enumitem}
	\usepackage{graphicx}
	\usepackage{tikz}
	\usepackage{algorithm}
	\usepackage{algorithmic}
	\usepackage{pgfplots}
	\usepackage{subfig}
	\usepackage[numbers]{natbib}
	\usepackage{todonotes}
	\usepackage{placeins}

	\font\msbm=msbm10

	\numberwithin{equation}{section}

	\theoremstyle{definition}
	\newtheorem{theorem}{Theorem}[section]
	\newtheorem{lemma}[theorem]{Lemma}

	\theoremstyle{definition}
	\newtheorem{definition}[theorem]{Definition}
	
	\newtheorem{remark}[theorem]{Remark}
	\def\mathbb#1{\hbox{\msbm{#1}}}

	\newtheorem*{remark*}{Remark}
	\newtheorem*{theorem*}{Theorem}

	\newcommand{\field}[1]{\ensuremath{\mathds{#1}}}

	\newcommand{\N}{\field N}

	\newcommand{\C}{\field{C}}

	\newcommand{\Z}{\field Z}
	
	\newcommand{\R}{\field R}
	\newcommand{\T}{\field T}

	\newcommand{\bproof}{\noindent {\bfseries Proof.}\, }
	\newcommand{\eproof}{\hspace*{\fill} \rule{3mm}{3mm}}

	\newcommand{\id}{\ensuremath{\mathrm{id}}}

	\newcommand{\rank}{\ensuremath{\mathrm{rank}}}

	\newcommand{\bh}{\mathbf{h}}
	\newcommand{\bk}{\mathbf{k}}
	\newcommand{\bj}{\mathbf{j}}
	\newcommand{\bx}{\mathbf{x}}
	
	\newcommand{\by}{\mathbf{y}}
	
	\newcommand{\bs}{\mathbf{s}}

	\newcommand{\bG}{\mathbf{G}}
	\newcommand{\bH}{\mathbf{H}}

	\newcommand{\bW}{\mathbf{W}}
	\newcommand{\bB}{\mathbf{B}}
	\newcommand{\bF}{\mathbf{F}}

	\newcommand{\cl}{\ensuremath{\mathcal L}}

	\newcommand{\be}{\begin{equation}}
	\newcommand{\ee}{\end{equation}}
	\newcommand{\beq}{\begin{eqnarray}}
	\newcommand{\eeq}{\end{eqnarray}}
	\newcommand{\beqq}{\begin{eqnarray*}}
	\newcommand{\eeqq}{\end{eqnarray*}}

	\newcommand{\cT}{\mathcal{T}}

	\newcommand{\Adeltaqn}{\ensuremath{A_{\Delta Q_n}}}
	\newcommand{\Binftyinftyrzero}{\ensuremath{\bB_{\infty, \infty}^{r_0} \torus}}
	\newcommand{\Binftyinftyr}{\ensuremath{\bB_{\infty, \infty}^{r} \torus}}
	\newcommand{\Bpastinftyr}{\ensuremath{\bB^r_{p^*, \infty} \torus}}
	\newcommand{\Bpinftyr}{\ensuremath{\bB^r_{p, \infty} \torus}}
	\newcommand{\Bppr}{\ensuremath{\bB_{p, p}^r \torus}}
	\newcommand{\Bpqr}{\ensuremath{\bB^r_{p, q} \torus}}
	
	\newcommand{\Btwotworone}{\ensuremath{\bB_{2, 2}^{r_1} \torus}}
	\newcommand{\Deltan}{\ensuremath{S_{\Delta Q_n}}}
	\newcommand{\Hpr}{\ensuremath{\bH^r_p \torus}}
	\newcommand{\Linfty}{\ensuremath{L_{\infty} \torus}}
	\newcommand{\Lp}{\ensuremath{L_p \torus}}
	\newcommand{\Ltwo}{\ensuremath{L_2 \torus}}
	
	\newcommand{\Sqn}{\ensuremath{S_{Q_n}}}
	\newcommand{\Wpr}{\ensuremath{\bW^r_p \torus}}
	
	\newcommand{\Wtwor}{\ensuremath{\bW^r_2} \torus}
	\newcommand{\Wtworone}{\ensuremath{\bW^{r_1}_2 \torus}}
	\newcommand{\alphvec}{\ensuremath{\boldsymbol{\alpha}}}
	\newcommand{\dm}[1]{\ensuremath{s_m\big(#1\big)}}
	
	\newcommand{\dmn}[1]{\ensuremath{s_{m_n}\big(#1\big)}}
	
	\newcommand{\hausdorff}{\ensuremath{\mathscr{A}}}
	\newcommand{\limJ}{\ensuremath{n_0}}
	
	\newcommand{\limL}{\ensuremath{n_1}}
	\newcommand{\torus}{\ensuremath{(\T^d)}}

\usepackage{xcolor}

\begin{document}

\title{Approximation of functions with small mixed smoothness in the uniform norm}
\author{Vladimir N. Temlyakov, Tino Ullrich\footnote{Corresponding author:
\url{tino.ullrich@mathematik.tu-chemnitz.de}} \\\\
University of South Carolina, Steklov Institute of Mathematics,\\ Lomonosov Moscow State University,\\ and Moscow Center for Fundamental and Applied Mathematics; \\
Faculty of Mathematics, 09107 Chemnitz, Germany}

\maketitle


\begin{abstract}
\sloppy In this paper we present results on asymptotic characteristics of multivariate function classes in the uniform norm. Our main interest is the approximation of functions with mixed smoothness parameter not larger than $1/2$.
Our focus will be on the behavior of the best $m$-term trigonometric
approximation as well as the decay of Kolmogorov and entropy numbers in the uniform norm. It turns out that these quantities share a few fundamental abstract properties like their behavior under real interpolation, such that they can be treated simultaneously. We start with proving estimates on finite rank convolution operators with range in a step hyperbolic cross. These results imply bounds for the corresponding function space embeddings by a well-known decomposition technique. The decay of Kolmogorov numbers have direct implications for the problem of sampling recovery in $L_2$ in situations where recent results in the literature are not applicable since the corresponding approximation numbers are not square summable.

\medskip

\small%

\end{abstract}


\section{Introduction} 
In this paper we provide new upper bounds for the best $m-$term trigonometric approximation ($\sigma_m$), the Kolmogorov numbers ($d_m$), and the entropy numbers ($e_m$) of multivariate function classes in the uniform norm. It is nowadays widely believed that the target space $\Linfty$ comes with additional difficulties and often requires new and involved techniques. Another challenge is the treatment of classes of periodic functions with small mixed smoothness (derivative or difference), where several questions concerning approximation and integration have not yet been settled. We make progress towards the solution of the Outstanding Open Problems 1.3, 1.6 and 1.7 in \cite{DuTeUl19}. The method used is rather general and reduces to a few common fundamental properties. This allows us to treat simultaneously the above asymptotic characteristics and make analogous statements for $s_m(T)$, where
$$
   ``s_m(T) \in \{d_m(T),e_m(T),\sigma_m(T)\}"
$$
and $T$ denotes an operator mapping into $L_{\infty}(\T^d)$. We follow the classical approach and start with new results for finite rank convolution operators $T = S_{Q_n}$, the orthogonal projection onto the trigonometric polynomials with frequencies in the dyadic step hyperbolic cross $Q_n \subset \Z^d$, defined by
\begin{align}\label{hyp1}
\varrho(\mathbf s) &:= \big\{ \mathbf k\in \Z^d:[ 2^{s_j-1}] \le
|k_j| < 2^{s_j}~,~ j=1,\dots,d \big\},\\
Q_n &:=   \bigcup\limits_{\mathbf \|\bs\|_1\le n}
\varrho(\mathbf s)\,.\label{hyp2}
\end{align}
Namely, for $2\leq p<\infty$ it holds
$$
	s_m(S_{Q_n}:\Lp \to \Linfty) \lesssim \left(\frac{2^n}{m}\right)^{\frac{1}{p}} n^{(d - 1)\left(1 - \frac{1}{p}\right) + \frac{1}{p}}\quad,\quad m \leq |Q_n|\,.
$$
The result is based on the common real interpolation properties of all three asymptotic characteristics in connection with a ``corner result'' due to Pajor, Tomczak-Jaegermann \cite{PaTo86}, Belinskii \cite{Bel98} and Dunker, K\"uhn, Linde, Lifshits \cite{DuLiKuLi98}. A corresponding corner result for the best $m$-term trigonometric approximation $\sigma_m$ in the univariate case was obtained by Belinskii \cite{Belinskii1998DecompositionTA}, who used a probabilistic technique, and in the multivariate case by Temlyakov \cite{Te15}, who used the greedy approximation technique.

It is well-known that the analysis of approximation problems for function classes  with small mixed smoothness involves several technical difficulties, see for instance \cite[Rem.\ 4.10]{Vy06} and \cite{MaUl20} for the study of entropy numbers. Similar difficulties have been already observed for the quantities of numerical integration, see \cite{UlUl16}, where the bounds look similar. Indeed, these quantities serve as lower bounds for Kolmogorov numbers in $L_\infty$,  which has been observed by Novak \cite{No86}. In this paper we give new asymptotic bounds for the classes $\Wpr$ and $\Hpr$, which are defined in Section 4 (see Definitions \ref{def2Sob}, \ref{defH} and Lemma \ref{FourierW}). In Theorems \ref{theorem:main1}, \ref{theorem:main2b} we obtain the results
$$
s_m(\operatorname{I}: \Wpr \to \Linfty) \lesssim m^{-r} (\log m)^{(d - 1) (1 - r) + r}
$$
as well as
$$
s_m(\operatorname{I}: \Hpr \to \Linfty) \lesssim m^{-r} (\log m)^{d-1 + r}\,
$$
in the ``small smoothness range'' $1/p<r<1/2$. In the endpoint situation $r=1/2$ we encounter an additional $(\log\log m)^{3/2}$ factor, see  Theorem \ref{theorem:main3} below. It is still open whether these bounds are sharp when $d>2$. The reader can find a brief discussion of the case $d=2$ in the Remark \ref{rem:d=2} below. Thus, we obtain new results on three asymptotic characteristics -- Kolmogorov numbers $d_m$, entropy numbers $e_m$, and best $m$-term approximations $\sigma_m$ -- for two kinds of classes $\bW^r_p$ and $\bH^r_p$ in the case of small smoothness $r\le 1/2$, when the error is evaluated in the uniform norm $L_\infty$. There is an extensive history of studying each of the above asymptotic characteristics. They were studied for large smoothness $r>1/2$, for classes
$\bW^r_p$, $\bH^r_p$, and for Besov classes $\bB^r_{p,\theta}$, where the error is evaluated in the $L_q$ norm, $1\le q\le \infty$. We refer the reader for a detailed historical discussion to the two recent books \cite{DuTeUl19} and \cite{VTbookMA}. For the $d_m$ see \cite[Sect.\ 4.3]{DuTeUl19}, and \cite[Sect.\ 5.3]{VTbookMA}. For the $e_m$ see \cite[Chapt.\ 6]{DuTeUl19} and \cite[Chapt.\ 7]{VTbookMA}. Finally, for the $\sigma_m$ see \cite[Chapt.\ 7]{DuTeUl19} and \cite[Chapt.\ 9]{VTbookMA}. In addition to the above books we mention the recent paper Romanyuk \cite{Rom19}.

We continue the investigation of asymptotic characteristics of classes of multivariate functions with small mixed smoothness started in \cite{TU20} on this topic. There we concentrated on the study of asymptotic characteristics from linear approximation theory -- the Kolmogorov widths. We pointed out some applications of new results on the Kolmogorov widths to the sampling recovery problem. In this paper focus is set on the study of asymptotic characteristics from nonlinear approximation theory -- sparse approximation with respect to the trigonometric system and entropy numbers. Those are interpreted as pseudo $s-$numbers sharing a few fundamental properties. We use a classical decomposition machinery (similar to the one used in \cite{TU20}), where we rely on finite rank operators ranging in subspaces of trigonometric polynomials with frequencies from hyperbolic crosses as building blocks. However, in contrast to \cite{TU20} we heavily apply well-known tools from interpolation theory of operators to analyze the finite rank operators. In \cite{TU20} an elementary approach is used to estimate {\em widths of function classes}, which is based on the application of a standard cutoff operator to dyadic building blocks. Certainly, deeply at the roots both approaches are related, since the cutoff operator is also used for computing $K$-functionals in real interpolation theory. However, technical realizations of these approaches are different and may be interesting for different communities.

Recent observations regarding the problem of optimal sampling recovery of functions in $L_2$ bring classes with small mixed smoothness to the focus again. Since several newly developed techniques only work for Hilbert-Schmidt operators \cite{KrUl19}, \cite{NaSchUl20} or, more generally, in situations where certain asymptotic characteristics (approximation numbers) are square summable \cite{KrUl20}, we need new techniques in situations where this is not the case. Especially in the range of small smoothness we are far away from square summability. Nevertheless, multivariate function classes of this type are of interest, since for instance a mixed H\"older-Zygmund regularity $r\leq 1/2$ falls into this scope. Recently, see \cite{temlyakov2020optimal}, the sampling recovery error in $L_2$ was directly related to the Kolmogorov numbers in $L_\infty$. It seems that, especially for the case of small smoothness, this represents the only available tool at the moment apart from sparse grid methods. Surprisingly, as an application of our results on Kolmogorov numbers we show that any sparse grid technique performs asymptotically worse by a $\log$-factor with exponent growing with the dimension $d$. This motivates further research in finding better constructive sampling algorithms.

The paper is organized as follows. In Sections 2 and 3 we define the asymptotic characteristics of interest in a framework of operators and pseudo $s$-numbers. This notion goes back to Pietsch \cite{Pietsch78}. We particularly pay attention to the real interpolation properties. Section 4 deals with the relevant function spaces with bounded mixed derivative or difference. Here we also give a new real interpolation formula. Afterwards in Section 5 we establish first results for the orthogonal projection operators with respect to the (trigonometric) step hyperbolic crosses. These estimates are used to obtain the main results in Section 6 for function space embeddings into $\Linfty$. Finally, in Section 7 we discuss the obtained results and give applications for the problem of sampling recovery.
\paragraph{Notation.} As usual $\N$ denotes the natural numbers, $\N_0:=\N\cup\{0\}$, $\Z$ denotes the integers, $\R$ the real numbers and \(\R_+\) the non-negative real numbers and $\C$ the complex numbers. $\C^n$ denotes the complex $n$-space. By $\T^d$ we denote the torus represented by the interval $[0,2\pi]^d$. Vectors or vector indices are usually typesetted boldface with, e.g., $\bx,\by\in \T^d$ or $\bs \in \N_0^d$. If not indicated otherwise  $\log(\cdot)$ denotes the natural logarithm of its argument. For $1\leq p\leq \infty$ and $\bx\in \C^n$ we denote $\|\bx\|_p := (\sum_{i=1}^n
|x_i|^p)^{1/p}$ with the usual modification in the case $p=\infty$ or $\bx$ being an infinite sequence. With $L_p(\T^d)$ we denote the space of all $p$-integrable $2\pi$-periodic complex-valued functions (equivalence classes) with $\int_{\T^d} |f(\bx)|^p\,d\bx<\infty$. For two sequences
$(a_n)_{n=1}^{\infty},(b_n)_{n=1}^{\infty}\subset \R_+$  we write $a_n \lesssim b_n$ if there exists a constant $c>0$,
such that $a_n \leq cb_n$ for all $n$. We will write $a_n \asymp b_n$ if $a_n \lesssim b_n$ and $b_n \lesssim a_n$.
If the constant $c$ depends on the dimension $d$ and smoothness $r$, we indicate it by $\lesssim_{r,d}$ and $\asymp_{r,d}$. For a linear operator between two normed spaces $X$, $Y$ we define its norm as
$\|T\|_{\mathcal{L}(X,Y)}:=\sup_{\|x\|_X \leq 1} \|Tx\|_Y$. The range of $T$ is defined as the subspace of $Y$ given by $\operatorname{range}(T) := T(X)$.


\section{Pseudo $s$-numbers}
\label{sec:kolmogorovNumbers}

In this section we introduce the asymptotic characteristics of interest, namely the Kolmogorov and entropy numbers as well as the error of best approximation with respect to an approximation scheme.

\begin{definition}[Kolmogorov numbers]
\label{definition:KolmogorovNumbers}
For Banach spaces $A, \, B$ and a linear operator $T: A \to B$, we define the
$m$-th Kolmogorov number as
$$
d_m(T: A \to B) := \inf \limits_{\substack{\dim L_m < m \\ L_m \subset B}} \sup
\limits_{\left\|a\right\|_{A} \leq 1} \inf \limits_{b \in L_m} \left\|T a
- b \right\|_B\quad,\quad m\in\N\,.
$$
\end{definition}
Let us start with the properties of the Kolmogorov numbers. The Kolmogorov numbers satisfy the following list of properties which make them a scale of $s$-numbers according to Pietsch \cite{Pie87}.
\begin{lemma}[Properties of Kolmogorov numbers]
\label{lemma:propertiesOfKolmogorovNumbers}
Let $A, \, B,\, C$ be Banach spaces and
$S, T \in \cl(A, B)$, $R \in \cl(B,C)$.
We have the following properties.
\begin{itemize}
\item[(K1)]	$\left\|T\right\|_{\cl(A, B)} = d_1(T) \geq d_2(T) \geq \dots \geq 0$,
\item[(K2)] For all $m_1, m_2 \in \N$, it holds
			$$
			d_{m_1 + m_2-1}(R \circ S) \leq d_{m_1}(R) d_{m_2}(S) \, .
			$$
\item[(K3)]	For all $m_1, m_2 \in \N$, it holds
			$$
			d_{m_1 + m_2-1}(S + T) \leq d_{m_1}(S) + d_{m_2}(T) \, .
			$$
\item[(K4)]	$d_m(T) = 0$ whenever $\rank{(T)} < m$.
\end{itemize}
\end{lemma}

Note that, except for (K4), these properties are shared by dyadic entropy numbers $(e_m)_m$ which we define below. To incorporate also dyadic entropy numbers into the framework Pietsch introduced the notion of {\em pseudo $s$-numbers}. We may use this notion here in a slightly different way.

\begin{definition}[Entropy numbers] Let $T:A\to B$ be a linear operator between two Banach spaces $A,B$. For $m\in \N$ the $m-$th dyadic entropy number of $T$ is defined as
$$
	e_m(T:A\to B) := \inf\Big\{\varepsilon>0~:~\exists b_1,...,b_{2^{m-1}}\in B \text{ such that } T(U_A) \subset \bigcup\limits_{k=1}^{2^{m-1}} (b_k + \varepsilon\cdot U_B) \Big\}\,,
$$
where $U_A$ and $U_B$ denote the unit balls in $A$ and $B$, respectively.
\end{definition}

Let us finally recall the definition of the asymptotic quantity measuring the best approximation with respect to an approximation scheme. This notion goes back to Pietsch \cite{Pie81} and includes the case of the best $m$-term approximation with respect to a dictionary $\mathcal{D}$. We will use it later for the multivariate trigonometric system. Let $X, Y$ denote arbitrary Banach spaces and let $(Y_n)_{n \in \N_0}$ denote a sequence of subsets of $Y$ satisfying
\begin{itemize}
	\item[(Y1)] $Y_0 = \{0\}$,
	\item[(Y2)] $Y_n \subset Y_{n+1}$, $n\in \N_0$,
	\item[(Y3)] $\lambda Y_n \subset Y_n$ for all $n\in \N_0$ and all scalars $\lambda$, and finally
	\item[(Y4)] $Y_n+Y_m \subset Y_{m+n}$\,.
\end{itemize}
\begin{definition}[Error of best approximation, \cite{Pie81}]
\label{defsigma} Let $X$ and $Y$ be as
above and let $(Y_n)_n$ denote a sequence of subsets in $Y$ fulfilling (Y1),...,(Y4) above. Let further $T:X\to Y$ denote a linear and bounded operator. Then we define the asymptotic characteristic
$$
	\sigma_m(T:X \to Y;(Y_n)_n):= \sup\limits_{\|x\|_X \leq 1}\inf\limits_{y \in Y_{m-1}}\|Tx-y\|_Y\quad,\quad m\in \N\,.
$$
\end{definition}

It turns out that counterparts of (K1), (K3) and (K4) hold true. (K2) has to be replaced by a weaker version (S2) which, however, is sufficient for our approach.

\begin{lemma}[Properties of $\sigma_m$]
Let $Z, \, X,\, Y$ be Banach spaces and
$S \in \cl(Z,X)$, $R, T \in \cl(X,Y)$. Let further $(Y_m)_m$ be a sequence of subsets in $Y$ fulfilling (Y1),...,(Y4) above. We have the following properties for $\sigma_m(T:X\to Y; (Y_k)_k)$.
\begin{itemize}
\item[(S1)]	$\left\|T\right\|_{\cl(X, Y)} = \sigma_1(T) \geq \sigma_2(T) \geq \dots \geq 0$,
\item[(S2)] For all $m \in \N$, it holds
			$$
			\sigma_{m}(T \circ S) \leq \sigma_{m}(R)\|S\|_{\mathcal{L}(Z,X)} \, .
			$$
\item[(S3)]	For all $m_1, m_2 \in \N$, it holds
			$$
			\sigma_{m_1 + m_2-1}(R + T) \leq \sigma_{m_1}(R) + \sigma_{m_2}(T) \, .
			$$
\item[(S4)]	If $\operatorname{range}(T) \subset Y_{m-1}$  then $\sigma_m(T) = 0$.
\end{itemize}
\end{lemma}
In the sequel we will often make statements for all three quantities at once. Then we will use the notation $s_m(T)$, where
$$
   ``s_m(T) \in \{d_m(T),e_m(T),\sigma_m(T)\}"\,.
$$
For technical reasons we put
\begin{equation}\label{s0}
	s_0(T):=0\,.
\end{equation}

\section{Real interpolation of pseudo $s$-numbers}
We first need the $K$-functional of a
Banach couple embedded into one joint Hausdorff space $\hausdorff$.
\begin{definition}[Peetre's $K$-functional, \cite{BeLo76}]
\label{definition:kFunctional}
For two Banach spaces $A_0, \, A_1$ which are jointly embedded into a common
Hausdorff space $\hausdorff$, we define for $a \in A_0 + A_1$
$$
K(t, a; A_0, A_1) = \inf \limits_{a = a_0 + a_1} \left(\left\|a_0\right\|_{A_0}
+ t \left\|a_1\right\|_{A_1}\right) \, .
$$
\end{definition}

The following interpolation results are well-known, see \cite[Sect.\ 11.6.8, 12.1.11]{Pietsch78}. Note, that a normed space $A$ is an {\it intermediate space} with respect to the couple $(A_0, A_1)$ if
$$
A_0 \cap A_1 \hookrightarrow A \hookrightarrow A_0+A_1\,,
$$
where ``$\hookrightarrow$'' indicates a continuous embedding. An intermediate space $A_{\theta}$ is of $K$-type $\theta$ if it satisfies
\be
\label{eqn:theta}
\sup \limits_{t > 0} t^{-\theta} K(t; a) \leq C \left\|a\right\|_{A_{\theta}}
\, .
\tag{$\Theta$}
\ee

\begin{theorem}[Interpolation of entropy and Kolmogorov numbers, \cite{Pietsch78}]
\label{theorem:interpolationOfKolmogorovNumbers}
Let $A_0$ and $A_1$ be embedded into the same Hausdorff space
$\hausdorff$. Let further $A_{\theta}$ be an intermediate space satsifying   condition \eqref{eqn:theta}.
Then, for any operator $T: A_0 + A_1 \to B$, one has
$$
e_{n + m - 1}(T: A_{\theta} \to B) \leq C \cdot e_n(T: A_0 \to B)^{1 - \theta}
e_m(T: A_1 \to B)^{\theta}
$$
and
$$
d_{n + m - 1}(T: A_{\theta} \to B) \leq C \cdot d_n(T: A_0 \to B)^{1 - \theta}
d_m(T: A_1 \to B)^{\theta} \, .
$$
\end{theorem}

The counterpart for the $(\sigma_m(T))_m$ numbers is straight-forward. Since we did not find such a result in the literature we decided to state it here explicitly and give a proof.
\begin{theorem}[Best approximation and interpolation]
\label{theorem:interpolationOfBestmterm}
Let $X_0, \, X_1$ and $X_{\theta}$ be embedded into the same Hausdorff space
$\hausdorff$. The intermediate space $X_{\theta}$ is supposed to satisfy \eqref{eqn:theta} with respect to the couple $(X_0,X_1)$. Then, we have for any linear operator $T: X_0 \to Y$, $T: X_1 \to Y$ with $Y$ and $(Y_k)_k$ as in Definition \ref{defsigma}
$$
\sigma_{n + m - 1}(T: X_{\theta} \to Y;(Y_k)_k) \leq C \cdot  ß\sigma_n(T: X_0 \to Y;(Y_k)_k)^{1 - \theta}
\sigma_m(T: X_1 \to Y;(Y_k)_k)^{\theta} \, .
$$

\bproof Let us abbreviate
$$
\sigma^0_n := \sigma_n(T: X_0 \to Y;(Y_k)_k), \qquad \sigma^1_m := (T: X_1 \to Y;(Y_k)_k)) \, .
$$
We clearly have for any $\varepsilon>0$, $x_0 \in X_0$ and $x_1 \in X_1$ elements $y_0 \in Y_{n-1}, y_1 \in Y_{m-1}$ such that
\be
\label{eqn:proofConditions}
\begin{aligned}
\left\|T x_0 - y_0\right\|_Y &\leq (1+\varepsilon)\sigma^0_n
\left\|x_0\right\|_{X_0} \, , \\
\left\|T x_1 - y_1\right\|_Y &\leq (1+\varepsilon)\sigma^1_m
\left\|x_1\right\|_{X_1} \, .
\end{aligned}
\ee
Let now $x \in X_{\theta}$
and $t > 0$. Then, for any $\delta>0$ there exist $x_0, \, x_1$ such that $x = x_0 + x_1$ and
$$
\left\|x_0\right\|_{X_0} + t \left\|x_1\right\|_{X_1} \leq C t^{\theta} \left\|x
\right\|_{X_{\theta}} (1 + \delta) \, .
$$
Put $t := \sigma^1_m/\sigma^0_n$ in the sequel (assuming $\sigma^0_n >0$, otherwise there is nothing to prove). Hence, due to \eqref{eqn:proofConditions}, there are $y_0, \, y_1$ such that
\begin{align*}
\left\|T x - (y_0 + y_1)\right\|_Y &\leq	\left\|T x_0 - y_0\right\|_Y +
											\left\|T x_1 - y_1\right\|_Y \\
									&\leq	(1 + \varepsilon) \left(\sigma^0_n \left\|
											x_0\right\|_{X_0} + \sigma^1_m \left\|x_1
											\right\|_{X_1}\right) \\
									&\leq	(1 + \varepsilon) \sigma^0_n \left(\left\|
											x_0\right\|_{X_0} + \frac{\sigma^1_m}{\sigma^0_n}
											\left\|x_1\right\|_{X_1}\right) \\
									&=		C (1 + \varepsilon)(1+\delta) \sigma^0_n \left(
											\frac{\sigma^1_m}{\sigma^0_n}\right)^{\theta} \left\|
											a\right\|_{X_{\theta}} \, .
\end{align*}
Put $y = y_0+y_1$ and observe by property (Y4) that $y \in Y_{m+n-2}$. Since $\varepsilon, \, \delta$ can be chosen arbitrarily small, we have
$$
\sigma_{n + m - 1}(T:X_{\theta} \to B) \leq C \sigma_n^{1 - \theta}(T:X_0 \to B) \cdot
\sigma_m^{\theta}(T:X_1 \to B) \, .
$$
\eproof
\end{theorem}


\section{Function spaces with small mixed smoothness}
For the subsequent definitions of the function classes of interest we refer to the recent books \cite{DuTeUl19} and \cite{VTbookMA} and to the references therein. We start with the definition of the univariate Bernoulli kernel. For fixed $\alpha \in \R$ we put
\[
F_{r,\alpha}(x):= 1+2\sum_{k=1}^\infty k^{-r}\cos (kx-\alpha \pi/2)\quad,\quad x\in \T\,.
\]
The corresponding multivariate Bernoulli kernels are defined via tensor products \index{Tensor product}. Let $\alphvec=(\alpha_1,\dots,\alpha_d) \in \R^d$ be fixed. Then
\begin{equation}\label{Bernoulli}
F_{r,\alphvec}(\bx):=\prod_{j=1}^dF_{r,\alpha_j}(x_j),\quad \bx=(x_1,\dots,x_d)\in
\T^d\,.
\end{equation}

Let us now proceed to classes with bounded mixed derivative. Note, that the parameter $\alphvec$ may be dropped (or set $\alphvec:=\mathbf{0}$) in the definition below in case $1<p<\infty$ since it leads to a family of equivalent norms, see also Lemma \ref{FourierW}. However, this is not the case in the endpoint cases $p=1$ and $p=\infty$.
\begin{definition} \label{def2Sob}
Let $r>0$, $\alphvec \in \R^d$ and $1\le p \le\infty$. Then $\bW^r_{p,\alphvec}$ is defined as the normed space of all $f\in L_p(\T^d)$ such that
\[
f=  F_{r,\alphvec}\ast
\varphi:=(2\pi)^{-d}\int_{\T^d}F_{r,\alphvec}(\bx-\by)\varphi(\by)d\by
\]
for some $\varphi \in L_p(\T^d)$, equipped with the norm
$\| \, f \, \|_{\bW^r_{p,\alphvec}}:= \|\varphi\|_p$.
\end{definition}

For the Littlewood-Paley characterization we need the building blocks $\delta_{\bs}(f,\bx)$, defined with \eqref{hyp1} by
\begin{equation}\label{defdelta}
   \delta_{\bs}(f,\bx):= \sum\limits_{\bk \in \varrho(\bs)} \hat{f}(\bk)e^{i\bk\cdot 	\bx}\,.
\end{equation}

\begin{lemma}\label{FourierW} If $1<p<\infty$ and $r>0$ then the norms $\|f\|_{\bW^r_{p,\alphvec}}$ with different $\alphvec$ are all equivalent to the Littlewood-Paley type norm
$$
\|f\|_{\Wpr} \asymp \Big\|\Big(\sum \limits_{\bs \in \N_0^d}
2^{r \|\bs\|_{1}} \left|\delta_{\bs}(f,\bx)\right|^2\Big)^{\frac%
{1}{2}}\Big\|_p\,.
$$
\end{lemma}
We now proceed with spaces with bounded mixed difference. Let $e$ be any subset of $\{1,...,d\}$. For multivariate functions $f:\T^d\to
\C$ and $\bh\in [0,1]^d$ the mixed first order difference operator $\Delta_{\bh}^{e}$
is defined by
\begin{equation*}
\Delta_{\bh}^{e} := \
\prod_{i \in e} \Delta_{h_i,i}\quad\mbox{and}\quad \Delta_{\bh}^{\emptyset} =  \operatorname{I},
\end{equation*}
where $\operatorname{I}f = f$ and $\Delta_{h_i,i}$ is the univariate first order
difference operator
$$
	\Delta_h g := g(\cdot+h)-g(\cdot)
$$
applied to the $i$-th variable of $f$ with the other variables kept
fixed. We first introduce spaces/classes $\bH^r_p$ of functions with bounded mixed difference. We restrict to first order difference operators since in this paper we are only interested in small smoothness.

\begin{definition}\label{defH}
Let $0<r<1$ and $1 \le p \le \infty$. We define the space $\bH^r_p$ as the set of all $f\in L_p(\T^d)$ such that for any $e \subset \{1,...,d\}$
\[
\big\|\Delta_{\bh}^{e}(f,\cdot)\big\|_p
\ \le \
C\, \prod_{i \in e} |h_i|^r
\]
for some positive constant $C$, and introduce the norm in this space
\[
\| \, f \, \|_{\bH^r_p} :=
\sum_{e \subset \{1,...,d\}}
\, | \, f \, |_{\bH^r_p(e)},
\]
where
\[
| \, f \, |_{\bH^r_p(e)} :=
\sup_{\bh \in (0,2\pi)^d} \,
\Big(\prod_{i \in e} |h_i|^{-r} \Big) \,\big\| \, \Delta_\bh^{e}(f,\cdot) \, \big\|_p\,.
\]
\end{definition}


For the purpose of the paper a characterization in terms of Fourier analytic building blocks is necessary. Since we also need to deal with $p=1$ and $p=\infty$ the blocks $\delta(f,\bx)$ will not be sufficient. We need the counterparts based on the classical de la Vall{\'e}e Poussin means, see \cite[Chapt.\ 2]{DuTeUl19}. Denote with $\mathcal{V}_m(t)$ the univariate de la Vall{\'e}e Poussin kernel
$$
	\mathcal{V}_m(t) = \frac{1}{m}\sum\limits_{k=m}^{2m-1} \mathcal{D}_k(t) = \frac{\sin(mt/2)\sin(3mt/2)}{m\sin^2(t/2)}\quad,\quad m\in \N\,.
$$
We further denote for $s \in \N_0$
$$
\mathcal{A}_s(t):=\left\{
	\begin{array}{rcl}\mathcal{V}_{2^s}(t)-\mathcal{V}_{2^{s-1}}(t)&:&s \geq 1,\\
					  \mathcal{V}_1(t)&:& s = 0\,.
					  \end{array}\right.
$$
In the multivariate case we use the tensorized version and define for $\bs \in \N_0^d$
$$
	\mathcal{A}_{\bs}(\bx):= \prod\limits_{i=1}^d
	\mathcal{A}_{s_i}(x_i)\quad, \quad \bx = (x_1,...,x_d)\,.
$$
Finally, the convolution operator $A_{\bs}(f,\cdot)$ is given by
\begin{equation}\label{Asdef}
	A_{\bs}(f,\cdot):=f\ast \mathcal{A}_{\bs}\,.
\end{equation}
For $1\leq p\leq \infty$ and $f\in L_p(\T^d)$ it holds $f = \sum_{\bs \in \N_0^d} A_{\bs}(f,\cdot)$ (with convergence in $L_p(\T^d)$) and, in particular,
\begin{equation}\label{Asbound}
	\|A_{\bs}:L_p(\T^d) \to L_p(\T^d)\| \asymp 1\quad,\quad \bs \in \N_0^d\,,
\end{equation}
see, e.g., \cite[(2.2.4)]{DuTeUl19}. For the following characterization we refer to
\cite[(3.3.2),(3.3.3)]{DuTeUl19}.

\begin{lemma}\label{Fourier1}
Let $0<r<1$. We have the following equivalent characterizations for $f\in \Lp$. \\
{\em (i)} If $1\leq p\leq \infty$ we have
\be
\begin{aligned}
\label{DefB}
\left\|f\right\|_{\Hpr} &\asymp \sup \limits_{\bs \in \N_0^d} \left\|A_\bs(f,\cdot)
\right\|_p 2^{r \|\bs \|_1} \, .
\end{aligned}
\ee
{\em (ii)} If $1 < p < \infty$ we have with \eqref{defdelta}
\be
\begin{aligned}
\label{DefB2}
\left\|f\right\|_{\Hpr} &\asymp \sup \limits_{\bs \in \N_0^d} \left\|\delta_\bs(f,\cdot)
\right\|_p 2^{r \|\bs \|_1} \, .
\end{aligned}
\ee
\end{lemma}

\begin{remark}\label{Fourier2} We will also need the refinement spaces $\Bpqr$, $1\leq q\leq \infty$, for technical reasons
\begin{equation}\label{f30}
 \|f\|_{\Bpqr} := \Big(\sum \limits_{\bs \in \N_0^d}
 2^{r \|\bs\|_1q}\left\| A_{\bs}(f,\cdot)\right\|^q_p \Big)^{\frac%
 {1}{q}} \, . \\
\end{equation}
In this notation we have $\bH^r_p(\T^d) = \Bpinftyr$ in the sense of equivalent norms. Note also, that in case $1<p<\infty$ we may replace $A_{\bs}(f,\cdot)$ by $\delta_{\bs}(f,\cdot)$ in \eqref{f30}. This together with Lemma \ref{FourierW} yields the identity $\bB^r_{2,2}(\T^d) = \bW^r_2(\T^d)$ in the sense of equivalent norms.
\end{remark}

Let us finally state a result on real interpolation of classes with bounded mixed difference which may be of interest on its own. Since the focus on the paper is on small smoothness we restrict here to smoothness parameters $r$ less than one. The theorem below also works for higher smoothness (using an isomorphism different from the Faber-Schauder system in the proof).

\begin{theorem}\label{Int_B} Let $2<p<\infty$, $0<r_0<1/2$ and $r_1 = r_0+1/2$. Then the following real interpolation formula
$$
	(\bB^{r_0}_{\infty,\infty}(\T^d), \bW^{r_1}_2(\T^d))_{\theta,p} = \bB^r_{p,p}(\T^d)
$$
(in the sense of equivalent norms) holds true if $\theta = 2/p$ and $r = r_0+1/p$.\\

\bproof First note that $\bW^{r_1}_2(\T^d) = \bB^{r_1}_{2,2}(\T^d)$ in the sense of equivalent norms. In \cite[Prop.\ 3.4 and 3.5]{HiMaOeUl16} it has been shown that for $1\leq p,q\leq \infty$ and $1/p<r<1+1/p$ it holds that
$$
	\|f\|_{B^r_{p,q}(\T^d)} \asymp \Big(\sum\limits_{\bj \in \N_{-1}^d}2^{\|\bj\|_1(r-1/p)q}
	\Big(\sum\limits_{\bk \in \mathbb{D}_{\bj}}|d^2_{\bj,\bk}(f)|^p\Big)^{q/p}\Big)^{1/q}\,.
$$
We used the notation from \cite[Sect.\ 3.3]{HiMaOeUl16}. Note, that $d^2_{\bj,\bk}(f)$ represents the Faber-Schauder coefficient of the corresponding $L_\infty$-normalized tensorized hat function with support $$[2^{-j_1}k_1,2^{-j_1}(k_1+1)]\times\cdots \times [2^{-j_d}k_d,2^{-j_d}(k_d+1)]$$ in the Faber-Schauder representation, see \cite[(3.5), (3.6)]{HiMaOeUl16}. As a direct consequence, the mapping $J_{r_0}$ (defined below) represents a common isomorphism between $\bB^{r_0+1/p}_{p,p}(\T^d)$ and $\ell_p$ for all $1\leq p\leq \infty$. Namely, we put
\begin{equation}
 \begin{split}
	J_{r_0} : \bB^{r_0+1/p}_{p,p}(\T^d) &\to \ell_p\\
	         f &\mapsto (2^{-r_0\|\bj\|_1}d^2_{\bj,\bk}(f))_{\bj,\bk}\,.
 \end{split}
\end{equation}
Hence, the real interpolation formula stated in the theorem is implied by the classical interpolation formula
$$
		(\ell_\infty,\ell_2)_{\theta,p} = \ell_p
$$
with $2<p<\infty$, $\theta = 2/p$, since $J_{r_0}$ maps all three occurring function spaces to either $\ell_\infty$, $\ell_2$ or $\ell_p$.
\eproof

\end{theorem}


\section{Convolution operators onto the step hyperbolic cross}
\label{Sect:Proj}
Let us refer to the definitions \eqref{defdelta} and \eqref{Asdef} of the dyadic block operators $\delta_{\bs} (f,\bx)$ and $A_{\bs}(f,\bx)$ based on the tensorized dyadic Dirichlet kernel and the tensorized de la Vall{\'e}e Poussin kernel, respectively. We further define for $n \in \N_0$ the hyperbolic cross operators
$$
\Sqn f:= \sum \limits_{\|\bs\|_1 \leq n} \delta_{\bs}(f,\cdot)
\qquad \text{and}
\qquad A_{Q_n} f = \sum \limits_{\|\bs\|_1 \leq n} A_{\bs}(f,\cdot).
$$
For the following facts we refer to \cite[Sect.\ 2.3]{DuTeUl19} and \cite[Sect.\  1.3.1, 3.2.5]{VTbookMA}. The range of $\Sqn$ represents the space of trigonometric polynomials with frequencies supported on the step hyperbolic cross $Q_n$, see \eqref{hyp1} and \eqref{hyp2}. We denote this space with
$$
	\cT(Q_n) := \operatorname{range}(\Sqn)\,.
$$
Due to the construction, $A_{Q_n}$ is invariant on $\cT(Q_n)$ such that we have
\begin{equation}\label{f2}
\Sqn = A_{Q_{n}}\circ \Sqn = \Sqn\circ A_{Q_n}\,.
\end{equation}
In contrast to $S_{Q_n}$, $A_{Q_n}$ ranges in a larger space of trigonometric polynomials. In particular, there is an integer $b\in \N$ such that
\begin{equation}\label{f3}
A_{Q_n} = S_{Q_{n+b}}\circ A_{Q_n} = A_{Q_n}
\circ S_{Q_{n+b}} \, .
\end{equation}
It is well-known
that in case $1<p<\infty$
\begin{equation}\label{normbound}
\|\Sqn:\Lp \to \Lp\| \asymp \|A_{Q_n}:\Lp \to \Lp\| \asymp 1 \quad,\quad
n\in \N_0 \, .
\end{equation}
Moreover, in case $p=\infty$ we have
\begin{equation}\label{Linfty_bound}
\|A_{Q_n}:\Linfty \to \Linfty \| \lesssim n^{d-1}\quad,\quad n\in \N \, ,
\end{equation}
see \cite[Lem.\ 4.2.3]{VTbookMA}.
We may also need the operators $\Deltan:=\Sqn -
S_{Q_{n-1}}$ and $\Adeltaqn:=A_{Q_n}-A_{Q_{n-1}}$ for $n\in \N$ (we set
$S_{\Delta Q_0} :=S_{Q_0}$ and $A_{\Delta Q_0} :=A_{Q_0}$).

In the sequel we are interested in the Kolmogorov, entropy numbers and best $m$-term trigonometric approximation of such a finite rank operator. More generally, we define for a finite set of frequencies $E \subset
\Z^d$ the corresponding projection operator $S_E: \Ltwo \to \Linfty$ given by
$$
S_E f(\bx) := \sum\limits_{\bk \in E}\hat{f}(\bk)e^{i\bk\cdot\bx }\,.
$$
In this context it is natural to study the best approximation error with respect to the multivariate trigonometric system with free spectrum, i.e., the best $m-$term trigonometric approximation of an operator $T$ defined by
\begin{equation}\label{sigma_m_trig}
	\sigma_m(T):=\sigma_m(T:\Lp \to \Linfty;(Y_n)_n)
\end{equation}
with
\begin{equation}\label{f31}
	Y_n := \Big\{t(\bx) = \sum\limits_{\bk \in \Lambda}c_{\bk}e^{i\bk\cdot \bx}~:~ |\Lambda| \leq n, c_{\bk} \in \C \Big\}\,.
\end{equation}

If not stated otherwise the quantity $\sigma_m$ will always be used in the context of best $m$-term trigonometric approximation in the sequel, see \eqref{sigma_m_trig}, \eqref{f31}. We may drop the $(Y_n)_n$ in the notation then.

In order to make use of the above interpolation results in Theorems \ref{theorem:interpolationOfKolmogorovNumbers} and \ref{theorem:interpolationOfBestmterm}, we may need an appropriate ``corner result''. The below bounds for entropy and Kolmogorov numbers are due to Pajor and Tomczak-Jaegermann \cite{PaTo86}, Belinskii \cite{Bel98}, see also 11.2.2 and 11.3.1 in \cite{TrBe04}, and Dunker, K\"uhn, Linde, Lifshits \cite{DuLiKuLi98}. The version for $(\sigma_m)_m$ is due to Temlyakov, see Theorem 2.6 in \cite{Te15}. For a univariate version of this result we refer to Belinskii \cite{Belinskii1998DecompositionTA}.

\begin{theorem}[\cite{PaTo86}, \cite{Bel98}, \cite{DuLiKuLi98}, \cite{Te15}]\label{TB} Let $E \subset \Z^d$ be a finite set such that
$$E \subset B_{\infty}(R) = \{\bx\in \R^d~:~\|\bx\|_\infty \leq R\}
$$
for some $R\geq 2$. Then we have

{\em (i)} for $``s_m \in \{d_m,\sigma_m\}"$
$$
		s_m(S_E: \Ltwo \to \Linfty) \lesssim \Big(\frac{|E|\log R}{m}\Big)^{1/2}\quad ,\quad m\leq |E|\,.
$$

{\em (ii)} For the entropy numbers it holds
$$
e_m(S_E: \Ltwo \to \Linfty) \lesssim
\left\{\begin{array}{rcl}
		\Big(\frac{|E|\log R}{m}\Big)^{1/2}&:&m \leq |E|,\\
		2^{-m/|E|}\sqrt{\log R}&:&m>|E|\,.
		\end{array}\right.
$$

\end{theorem}


The following bounds are direct consequences of Theorem \ref{TB} in connection with Theorem \ref{theorem:interpolationOfKolmogorovNumbers}. Note, that (i) below for $s_m = e_m$ is already known. It was obtained in \cite{Te20_5} to prove the Marcinkiewicz
type discretization theorems for the hyperbolic cross polynomials. The proof there is based on a different technique.

\begin{theorem}
\label{theorem:stepHyperbolicCrosses}
Let $2 \leq p < \infty$. Then, it holds for $1\leq m\leq |Q_n|$ and $``s_m \in \{d_m,e_m,\sigma_m\}"$
\begin{enumerate}
\item	$$
		s_m(\Sqn: \Lp \to \Linfty)
		\lesssim \left(\frac{2^n}{m}\right)^{\frac{1}{p}} n^{(d - 1)
		\left(1 - \frac{1}{p}\right) + \frac{1}{p}} \, .
		$$
\item If $r\geq 0$ then $$
		s_m(\Deltan: \Wpr \to \Linfty)
		\lesssim 2^{-rn}\left(\frac{2^n}{m}\right)^{\frac{1}{p}} n^{(d - 1)
		\left(1 - \frac{1}{p}\right) + \frac{1}{p}},
		$$
\item	and if $r>1/p$
		$$
		s_m(\Adeltaqn: \Bpinftyr \to \Linfty) \lesssim 2^{-r
		n} \left(\frac{2^n}{m}\right)^{\frac{1}{p}} n^{d - 1 + \frac{1}{p}} \, .
		$$
\end{enumerate}
\bproof
Let us start proving (i). It is well-known for the real interpolation method
$(\cdot, \cdot)_{\theta,q}$ that
$\left(\Linfty, \Ltwo\right)_{\theta, p} = \Lp$ whenever
$$
\frac{1}{p} = \frac{1 - \theta}{\infty} + \frac{\theta}{2},
$$
which means $\theta = 2/p$. Hence, we have with $A_{\theta} = \Lp$ the condition \eqref{eqn:theta} fulfilled.
So, we may interpolate the numbers $s_m$ according to Theorems
\ref{theorem:interpolationOfKolmogorovNumbers} and \ref{theorem:interpolationOfBestmterm}. This gives for the operator $A_{Q_n}$
\begin{align*}
&\dm{A_{Q_n}: \Lp \to \Linfty} \\
&~~~~~~~~~\lesssim \left\|A_{Q_n}: \Linfty \to \Linfty
\right\|^{1 - \theta} \cdot \dm{A_{Q_n}: \Ltwo \to \Linfty}^{\theta}\\
&~~~~~~~~~\lesssim n^{(d-1)(1-\theta)}\dm{A_{Q_n}: \Ltwo \to \Linfty}^{\theta} \, ,
\end{align*}
where we used \eqref{Linfty_bound}. We continue with the first identity in
\eqref{f3}, where $S_{n+b}: \Ltwo \to \Linfty$ and $A_{Q_n}: \Ltwo \to
\Ltwo$. Then, (K2) and \eqref{normbound} yield
$$
\dm{A_{Q_n}: \Ltwo \to \Linfty}^{\theta} \lesssim \dm{S_{Q_{n+b}}: \Ltwo \to
\Linfty}^{\theta} \, .
$$
Applying Theorem \ref{TB} to the right-hand side together with $\theta = 2/p$ yields for $m \leq |Q_n|$
\begin{equation}\label{dm_A}
\dm{A_{Q_n}: \Lp \to \Linfty} \lesssim
\Big(\frac{2^n}{m}\Big)^{\frac{1}{p}}n^{(d-1)(1-\frac{1}{p})+\frac{1}{p}} \, .
\end{equation}
Using the first identity in \eqref{f2} together with the properties (K2), (S2) and  \eqref{normbound} gives
$$
\dm{\Sqn:\Lp \to \Linfty} \lesssim \dm{A_{Q_{n}}:\Lp \to
\Linfty} \,,
$$
where the right-hand side can be bounded by \eqref{dm_A}. This proves (i).

For (ii) observe that by Lemma \ref{FourierW}
$$
\dm{\Deltan:\Wpr \to \Linfty} \asymp 2^{-rn}
\dm{\Deltan:\Lp \to \Linfty},
$$
which will be bounded using (i).

As for (iii), we have by the real interpolation formula in Theorem \ref{Int_B} with
$r_0 = r-1/p$, $r_1 = r-1/p+1/2$ and $\theta = 2/p$
$$
\left(\Binftyinftyrzero, \Btwotworone\right)_{\theta,
p} = \Bppr \quad , \quad 2 < p < \infty \,.
$$
As a consequence, we obtain the condition \eqref{eqn:theta} for $A_{\theta} = \Bppr$ with respect to the above couple. Interpolating Kolmogorov numbers according to
\Cref{theorem:interpolationOfKolmogorovNumbers} gives
\begin{align}\label{f9}
\dm{\Adeltaqn: \Bppr \to \Linfty} \lesssim& \|\Adeltaqn: \Binftyinftyrzero
\to \Linfty\|^{1 - \theta} \\ & ~\cdot\nonumber
\dm{\Adeltaqn: \Btwotworone \to \Linfty}^{\theta}.
\end{align}
By Lemma \ref{Fourier1} together with \eqref{Asbound} and $\theta = 2/p$ we find
\begin{equation}\label{f10}
\|\Adeltaqn:\Binftyinftyrzero
\to \Linfty\|^{1 - \theta} \lesssim n^{(d - 1) \left(1 - \frac{2}{p}\right)}
2^{-r_0 n \left(1 - \frac{2}{p}\right)}.
\end{equation}
Since $\Btwotworone = \Wtworone$ in the sense of equivalent norms, we
may use \eqref{dm_A} and plug the result together with \eqref{f10} into \eqref{f9}. This
yields
\begin{align}
\nonumber
\dm{\Adeltaqn: \Bppr \to \Linfty} &\lesssim
n^{(d - 1) \left(1 - \frac{2}{p}\right)}
2^{-r_0 n\left(1 - \frac{2}{p}\right)}2^{-r_1 n\frac{2}{p}} \left[\frac{2^n
\cdot n^{d - 1}n}{m}\right]^{\frac{1}{p}}\\
&\asymp 2^{-rn}\Big(\frac{2^n}{m}\Big)^{\frac{1}{p}}n^{(d-1)(1-\frac{1}{p})+
\frac{1}{p}} \, .\label{f11}
\end{align}
Finally, by Lemma \ref{Fourier1}, (ii), Remark \ref{Fourier2} and \eqref{Asbound},
we see that $\|S_{Q_{n+b}}-S_{Q_{n-b}}:\Bpinftyr \to \Bppr\| \asymp n^{\frac{d-1}{p}}$. Using $\Adeltaqn = \Adeltaqn\circ (S_{Q_{n+b}}-S_{Q_{n-b}})$
together with (K2), (S2) gives
\begin{align*}
\dm{\Adeltaqn:\Bpinftyr \to \Linfty} &\lesssim n^{\frac{d-1}{p}}
\dm{\Adeltaqn:\Bppr \to \Linfty}\\
&\lesssim 2^{-rm}\Big(\frac{2^n}{m}\Big)^{\frac{1}{p}}n^{d-1 + \frac{1}{p}} \, ,
\end{align*}
where we used \eqref{f11} in the last step. \eproof
\end{theorem}


\section{Embeddings into $\Linfty$}
\label{main}
Let us present here our main results for embeddings of Sobolev and
H\"older-Nikolskii spaces with small mixed smoothness into $\Linfty$.

\begin{theorem}
\label{theorem:main1}
Let $2 < p < \infty, \, \frac{1}{p} < r < \frac{1}{2}$. Then, for $``s_m \in \{d_m,e_m,\sigma_m\}"$ we have
$$
\dm{\operatorname{I}: \Wpr \to \Linfty} \lesssim m^{-r} (\log m)^{(d - 1) (1 - r) + r} \, .
$$

\bproof
We decompose the identity operator
$$
\operatorname{I} = \sum \limits_{n = 0}^{\infty} \Deltan \, ,
$$
where the $\Deltan$ are the operators defined above. Using (K3), (S3) we
have that
\begin{equation}\label{f101}
s_m(\operatorname{I}) \leq \sum \limits_{n = 0}^{\infty} \dmn{\Deltan}
\end{equation}
with $m = \sum_{n = 0}^{\infty} m_n$. Note that this implies $m_n = 0$ for $n\geq n_0$ since $m_n \in \N_0$. Hence, in view of \eqref{s0}, the above sum makes sense. We further decompose into three parts
\begin{equation}\label{f20}
s_m(\operatorname{I}) \leq \sum \limits_{n = 0}^{\limJ} \dmn{\Deltan} + \sum \limits_{n = \limJ}^{\limL} \dmn
{\Deltan} + \sum \limits_{n = \limL}^{\infty} \dmn{\Deltan} \, .
\end{equation}
Let us consider the first sum in \eqref{f20}. The following argument only works for $``s_m \in \{d_m,\sigma_m\}"$ since a counterpart of (K4) or (S4) is not available for entropy numbers. We will indicate the necessary modification for $s_m = e_m$ below. Let $\limJ$ be the largest number such that
$$
\sum \limits_{n = 0}^{\limJ} \rank\big(\Deltan\big) \leq m
$$
and put $m_n := \rank\big(\Deltan\big) + 1$. Due to property (K4) and (S4) in  Lemma \ref{lemma:propertiesOfKolmogorovNumbers}, we make the first sum disappear.
As for the second sum, we choose
$$
m_n := \lfloor2^n 2^{(\limL - n) \kappa} \limL^{-(d - 2)}\rfloor
$$
and $\limL$ such that
\be
\label{eqn:main1}
\frac{2^{\limL}}{\limL^{d - 2}} \asymp m \, .
\ee
Clearly,
$$
\sum \limits_{n = \limJ}^{\limL} 2^n 2^{(\limL - n) \kappa} \limL^{-(d - 2)} \asymp \frac{2^{\limL}}{\limL^
{d - 2}} \asymp m \, .
$$
Here, $\kappa$ is chosen such that $2 r < \kappa < 1$. Let us decompose as follows
$$
\Deltan = \Deltan \circ I
$$
with $\operatorname{I}:{\Wpr \to \Wtwor}$ and get
\begin{align*}
\dmn{\Deltan: \Wpr \to \Linfty}	&\lesssim	\left\|\operatorname{I}: \Wpr \to \Wtwor \right\| \\
								&\;			\cdot \; \, \,\dmn{\Deltan: \Wtwor \to
											\Linfty} \\
								&\lesssim	2^{-r n} \big(2^{-(\limL - n) \kappa}
											\limL^{d - 2}\big)^{\frac{1}{2}} n^{\frac%
											{d - 1}{2} + \frac{1}{2}} \, ,
\end{align*}
where we used \Cref{theorem:stepHyperbolicCrosses} (actually
\cite[Theorem~11.3.1]{TrBe04}). Summing up with $\limJ \leq n \leq \limL$ gives
\be
\label{eqn:main2}
\sum \limits_{n = \limJ}^{\limL} \dmn{\Deltan: \Wpr \to \Linfty} \lesssim 2^{-r \limL} \limL^{\frac{d
- 2}{2}} \limL^{\frac{d - 1}{2} + \frac{1}{2}} \, .
\ee
Using the fact that $\frac{2^{\limL}}{\limL^{d - 2}} \asymp m$, we have
\be
\begin{aligned}
\label{eqn:main3}
\eqref{eqn:main2}
&\lesssim m^{-r} \limL^{-(d - 2) r + \frac{d - 2}{2} + \frac{d - 1}{2} + \frac{1}{2}} \\
&\asymp m^{-r} \limL^{(d - 1) (1 - r) + r} \, .
\end{aligned}
\ee
Now we care for the third sum and choose
$$
m_n := \lfloor m \cdot 2^{(\limL - n) \zeta} \rfloor\, ,
$$
where $\zeta>0$ is chosen such that
$$
\frac{\zeta}{p} < r - \frac{1}{p} \, .
$$
Clearly, $\sum_{n = \limL}^{\infty} m_n \asymp m$. By the results from
\Cref{theorem:stepHyperbolicCrosses}, we obtain
$$
\dmn{\Deltan: \Wpr \to \Linfty} \lesssim \left[\frac{2^n}{m} \cdot 2^{(n - \limL)
\zeta}\right]^{\frac{1}{p}} 2^{-r n} n^{(d - 1) \left(1 - \frac{1}{p}\right) +
\frac{1}{p}} \, .
$$
Summing over $n$ in the range $n = \limL, \limL + 1, \dots$ gives (taking \eqref{s0} into account)
\be
\label{eqn:mainSumOverN}
\sum \limits_{n = \limL}^{\infty} \dmn{\Deltan} \lesssim \left(\frac{2^{\limL}}{m}\right)^
{\frac{1}{p}} 2^{-r \limL} \cdot \limL^{(d - 1) \left(1 - \frac{1}{p}\right) + \frac{1}{p}}
\, .
\ee
Because of \eqref{eqn:main1}, we have
\begin{align*}
\eqref{eqn:mainSumOverN}
&\lesssim \limL^{\frac{d - 2}{p}} 2^{-r \limL} \limL^{(d - 1) \left(1 - \frac{1}{p}\right)
+ \frac{1}{p}} \\
&\asymp 2^{-r \limL} \limL^{d - 1} \\
&\asymp m^{-r} \limL^{d - 1} \limL^{-(d - 2) r} \\
&\asymp m^{-r} \big(\log m\big)^{(d - 1) (1 - r) + r} \, .
\end{align*}
This, combined with \eqref{eqn:main3}, gives the result of the theorem for $``s_m \in \{d_m, \sigma_m\}"$.

We finally comment on the estimate of the first sum in \eqref{f20} in case of entropy numbers. We modify the argument as follows: Instead of choosing $m_n = \rank(S_{\Delta Q_n})+1$ we choose
\begin{equation}\label{f21}
	m_n := \lfloor \rank(S_{\Delta Q_n})2^{(\limJ - n)\varepsilon}\rfloor\quad,\quad n = 1,...,\limJ\,,
\end{equation}
with $0<\varepsilon < 1$. This gives
$$
	\sum\limits_{n=0}^{\limJ} m_n \asymp \rank(S_{\Delta Q_{\limJ}}) \asymp 2^{\limJ} \limJ^{d-1}\,,
$$
where we choose $\limJ$ such that $2^{\limJ} \limJ^{d-1} \asymp m$. By (K2) and Theorem \ref{TB} we obtain (note that $m_n>\rank(S_{\Delta Q_n})$)
$$
	e_{m_n}(S_{\Delta Q_n}:\Wpr \to \Linfty)\lesssim
	2^{-rn}2^{-2^{(\limJ - n)\varepsilon}}n^{1/2}\,.
$$
Summing over $n = 0,...,\limJ$ yields
\begin{align*}
	&\sum\limits_{n=0}^{n_0} e_{m_n}(S_{\Delta Q_n}:\Wpr \to \Linfty)
	\\&~~~~~~~~\lesssim 2^{-r \limJ} \limJ^{1/2}\asymp m^{-r}(\log m)^{r(d-1)+1/2} \lesssim m^{-r}(\log m)^{(d-1)(1-r)+r}\,.
\end{align*}
This finishes the proof.
\eproof
\end{theorem}

\begin{theorem}
\label{theorem:main2b}
Let $2 < p \leq \infty$ and $\frac{1}{p} < r < \frac{1}{2}$. Then, for $``s_m \in \{d_m,e_m,\sigma_m\}"$
$$
\dm{\operatorname{I}: \Hpr \to \Linfty} \lesssim m^{-r} \big(\log m\big)^{d - 1 + r} \, .
$$

\bproof This time, we decompose the identity using the operators $\Adeltaqn$. Let us first deal with the case $p<\infty$. Applying again property (K3) and (S3) we find (taking \eqref{s0} into account)
\be
\label{eqn:main5}
\dm{\operatorname{I}: \Bpinftyr \to \Linfty} \lesssim \sum \limits_{n = 0}^{\limJ} \dmn{\Adeltaqn} +
\sum \limits_{n = \limJ}^{\limL} \dmn{\Adeltaqn} + \sum \limits_{n = \limL}^{\infty} \dmn
{\Adeltaqn} \, .
\ee
We argue analogously as in the proof of Theorem \ref{theorem:main1} for the first sum.
For the second sum we choose
\begin{equation}\label{f41}
m_n = \lfloor 2^n 2^{(\limL - n) \kappa} \limL \rfloor
\end{equation}
with $2 r < \kappa < 1$ and $n_1$ such that $2^{\limL} \limL \asymp m$. Hence,
\be
\label{eqn:main4}
\sum \limits_{n = \limJ}^{\limL} m_n \asymp 2^{\limL} \limL \asymp m \, .
\ee
Then we decompose
\begin{equation}\label{f40}
	A_{\Delta Q_n} = S_{Q_{n+b}}\circ A_{\Delta Q_n}\,.
\end{equation}
This gives
\begin{align*}
\dmn{\Adeltaqn: \Bpinftyr \to \Linfty}
&\lesssim \left\|A_{\Delta Q_n}: \Bpinftyr \to \bB^r_{2,2}(\T^d) \right\| \\
&~~~~\cdot s_{m_n}(S_{Q_{n+b}}: \bB^r_{2,2}(\T^d) \to \Linfty) \\
&\lesssim 2^{-r n} \left(\frac{2^n}{m_n}\right)^{\frac{1}{2}} n^{d - 1 + \frac%
{1}{2}} \, ,
\end{align*}
where we used \Cref{theorem:stepHyperbolicCrosses}, (ii) for estimating $s_{m_n}$. To estimate $\|A_{\Delta Q_n}\|$ we used Lemma \ref{Fourier1} together with \eqref{Asbound}. Inserting \eqref{f41} yields
$$
\dmn{\Adeltaqn: \Bpinftyr \to \Linfty} \lesssim 2^{-r n} \left[2^{(n - \limL) \kappa}
\cdot \limL^{-1}\right]^{\frac{1}{2}} n^{d - 1 + \frac{1}{2}} \, .
$$
Summation over $n = \limJ, \dots , \limL$ leads to
$$
\sum \limits_{n = \limJ}^{\limL} \dmn{\Adeltaqn} \lesssim \frac{2^{-r \limL}}{\sqrt{\limL}} \limL^{d
- 1 + \frac{1}{2}} \, .
$$
Since $m \asymp 2^{\limL} \limL$, due to \eqref{eqn:main4}, we get
\be
\label{eqn:main*}
\sum \limits_{n = \limJ}^{\limL} \dmn{\Adeltaqn} \lesssim 2^{-r \limL} \limL^{d - 1 + r} \asymp
m^{-r} \big(\log m\big)^{d - 1 + r} \, .
\ee

We finally deal with the last sum in \eqref{eqn:main5}. Indeed, by choosing
$$
m_n := \lfloor m \cdot 2^{(\limL - n) \zeta} \rfloor
$$
with $r - \frac{1}{p} - \frac{\zeta}{p} > 0$, we have by
\Cref{theorem:stepHyperbolicCrosses} (for $m_n>0$)
\begin{align*}
\dmn{\Adeltaqn: \Bpinftyr \to \Linfty}
&\lesssim 2^{-r n} \left(\frac{2^n}{m_n}\right)^{\frac{1}{p}} n^{d - 1 + \frac%
{1}{p}} \\
&\asymp \left(\frac{2^n}{m} 2^{(n - \limL) \zeta}\right)^{\frac{1}{p}} 2^{-r n}
n^{d - 1 + \frac{1}{p}} \, .
\end{align*}
Summing up over $n = \limL, \limL + 1, \dots$ yields (taking \eqref{s0} into account)
$$
\sum \limits_{n = \limL}^{\infty} \dmn{\Adeltaqn} \lesssim 2^{-r \limL} \cdot \left(
\frac{2^{\limL}}{m}\right)^{\frac{1}{p}} \limL^{d - 1 + \frac{1}{p}} \,.
$$
By \eqref{eqn:main4} we get $\frac{2^{\limL}}{m} \asymp \frac{1}{\limL}$. Hence, we obtain
$$
\sum \limits_{n = \limL}^{\infty} \dmn{\Adeltaqn} \lesssim m^{-r} \big(\log m\big)^
{d - 1 + r} \, .
$$
Together with \eqref{eqn:main*}, this proves the theorem in case $p < \infty$.
For the case $p = \infty$ we use the bounded embedding $\operatorname{I}:{\Binftyinftyr \to
\Bpastinftyr}$, where $p^* < \infty$ is chosen such that $r > \frac{1}{p^*}$.
This gives
\begin{align*}
\dm{\operatorname{I}: \Binftyinftyr \to \Linfty}
&\leq \left\|\operatorname{I}: \Binftyinftyr \to \Bpastinftyr\right\| \\
&~~~~~\cdot \dm{\operatorname{I}: \Bpastinftyr \to \Linfty} \\
&\lesssim m^{-r} \big(\log m\big)^{d - 1 + r} \, ,
\end{align*}
where we used the result for $p^* < \infty$.

Again, we comment on the necessary modifications in case of $s_m = e_m$. Let us consider the first sum in \eqref{eqn:main5} again and use \eqref{f40}\,. We choose  $m_n$ and $n_0$ as after \eqref{f21}. By the counterpart of (K2) for entropy numbers we find
\begin{align*}
	&e_{m_n}(A_{\Delta {Q_n}}:\Bpinftyr \to \Linfty) \\
	&~~~~~~\leq \|A_{\Delta Q_n}:\Bpinftyr \to \bB^r_{2,2}(\T^d)\|e_m(S_{Q_{n+b}}:\bB^r_{2,2}(\T^d) \to \Linfty)\\
	&~~~~~~\leq n^{(d-1)/2}2^{-rn}2^{-2^{(\limJ-n)\varepsilon}}n^{1/2}\,,
\end{align*}
where we applied Theorem \ref{TB}, (ii). Summing over $n=0,...,\limJ$ yields
$$
	\sum\limits_{n=0}^{n_0}e_{m_n}(A_{\Delta Q_n})\lesssim 2^{-r \limJ} \limJ^{(d-1)/2+1/2} \asymp m^{-r}(\log m)^{(d-1)(r+1/2)+1/2}\lesssim
	m^{-r}(\log m)^{d-1+r}\,.
$$
This concludes the proof.
\eproof
\end{theorem}

For the endpoint situation $r=1/2$ we obtain an additional $(\log\log m)^{3/2}$ factor in the upper bounds.

\begin{theorem}[Endpoint cases]\label{theorem:main3} Let $``s_m \in \{d_m,e_m,\sigma_m\}"$.\\ {\em (i)} If $2<p<\infty$ and $r=1/2$ then
$$
   \dm{\operatorname{I}: \Wpr \to \Linfty} \lesssim m^{-r} (\log m)^{(d - 1) (1 - r) + r}(\log\log m)^{r+1} \, .
$$
{\em (ii)} If $2<p\leq \infty$ and $r=1/2$ then
$$
\dm{\operatorname{I}: \Hpr \to \Linfty} \lesssim m^{-r} \big(\log m\big)^{d - 1 + r}(\log\log m)^{r+1} \, .
$$
\bproof We use the same decomposition of the identity operator as above. The first and third sum will be treated analogously. In the second sum it is not possible to choose $\kappa < 1$. We choose $\kappa = 1$ but pay a $\log(\limL)$ in both summations. Rephrasing the final bound in terms of $m$ yields an additional $(\log \log m)^{r+1}$ factor.  \eproof
\end{theorem}

\begin{remark}[$d=2$]\label{rem:d=2} {\em (i)} We would like to emphasize that in  Theorem \ref{theorem:main2b}, when $d=2$, we actually do not need the middle sum ranging over $[n_0,n_1]$ in \eqref{eqn:main5}. Hence, the restriction $r\leq 1/2$ does not play a role here. This results in
\begin{equation}\label{d=2}
 s_m(\operatorname{I}: \bH^r_{p}(\T^2) \to L_{\infty}(\T^2)) \lesssim m^{-r} \big(\log m\big)^{1 + r}
\end{equation}
for all $r>1/p$ and $2\leq p\leq \infty$. Compared to Theorem \ref{theorem:main3}, (ii) we do not have a $\log\log$-term here for $r=1/2$. In addition, together with Theorem 7.8.4 from \cite{VTbookMA} (see also \cite[Thm.\ 6.3.4]{DuTeUl19} and the references therein) and Carl's inequality \cite{Ca81} we get the correct order in case $d=2$ for Kolmogorov and entropy numbers. Namely for $2 \leq p\leq \infty$ and $r>1/p$ it holds
$$
	 e_m(\operatorname{I}: \bH^r_{p}(\T^2) \to L_{\infty}(\T^2))
	 \asymp d_m(\operatorname{I}: \bH^r_{p}(\T^2) \to L_{\infty}(\T^2))
	 \asymp m^{-r} \big(\log m\big)^{1 + r} \, .
$$
The result for entropy numbers is true for $1 \leq p \leq \infty$, $r >1/p$, see \cite[Thm.\ 7.8.4]{VTbookMA}.
Note, that for the $\bW^r_p$ classes the correct order of decay for $d=2$ of the Kolmogorov and entropy numbers is only known in case of large smoothness $r>1/2$, see \cite[Chapt.\ 11]{TrBe04} and also \cite{DuLiKuLi98}.

{\em (ii)} The upper bound in \eqref{d=2} also includes the error of best $m$-term trigonometric approximation $\sigma_m$. Together with the lower bounds from \cite[Thm.\ 3.3]{Te15} we have in case $d=2$, $2\leq p \leq \infty$ and $r>1/p$
\begin{equation*}
 m^{-r}\log m \lesssim \sigma_m(\operatorname{I}: \bH^r_{p}(\T^2) \to L_{\infty}(\T^2)) \lesssim m^{-r} \big(\log m\big)^{1 + r}\,.
\end{equation*}

\end{remark}

\section{Applications and discussion}

In this section we comment on applications of the above results and add a discussion on possible future research and open problems motivated by our considerations. We can say right here that we already have made progress on the Outstanding Open Problems in \cite{DuTeUl19}, especially 1.3, 1.6, 1.7. In addition, we discuss consequences for sampling recovery in $L_2$. Furthermore, we comment on the use of finite dimensional subspaces generated by hyperbolic wavelets as buidling blocks and wavelet type dictionaries for best $m$-term approximation.

\paragraph{Entropy and Kolmogorov numbers.} Entropy numbers for mixed smoothness embeddings have been investigated by several authors in the literature, see \cite[Chapt.\, 6]{DuTeUl19}. Among many others, Vyb{\'i}ral \cite{Vy06} investigated the behavior of entropy numbers in $\bB^s_{p,q}$-spaces, see Remark \ref{Fourier2}, using wavelet building blocks. In addition, the authors in \cite{DuLiKuLi98} managed to prove a counterpart of the corner result in Theorem \ref{theorem:stepHyperbolicCrosses}, (i), for $p=2$, $s_m = d_m$ and $S_{Q_n}$ replaced by the corresponding hyperbolic Haar wavelet projection.

Let us comment on this technique here and how it can be applied for the uniform
norm estimates. Technically, instead of trigonometric polynomials one may also use a univariate wavelet system $\{\psi_I = \psi((\cdot-x_I)/|I|)~:~I \in \mathcal{I}, |I| \leq 1\}$, where $\mathcal{I}$ is the set of dyadic intervals $I$ with midpoints $x_{I} = k2^{-j}$, $k, j \in \Z$. We further consider the corresponding multivariate (tensorized) system
\begin{equation}\label{f50}
  \mathcal{D} = \Big\{\psi_{\mathbf{I}}(\bx) = \prod\limits_{j=1}^d \psi_{I_j}(x_j) ~:~\mathbf{I} = I_1\times\cdots \times I_d, |I_j| \leq 1, j=1,...,d\Big\}
\end{equation}
and define the orthogonal projection on the hyperbolic layers
$$
	\tilde{S}_n f := \sum\limits_{\substack{\mathbf{I} \subset [0,1]^d\\ |\mathbf{I}| = 2^{-n} }} \langle f, \psi_{\mathbf{I}} \rangle \psi_{\mathbf{I}}\,.
$$
This operator replaces the above $S_{\Delta Q_n}$. Then we have the decomposition of the identity operator $\operatorname{I} = \sum_{n=0}^{\infty} \tilde{S}_n\,.$
We assume, that the wavelet system is sufficiently smooth, compactly supported and has good decay properties. We also need the finite-dimensional block result
\begin{equation}\label{f39}
	s_m(\id:\ell^N_p \to \ell^N_\infty) \asymp \Big[\frac{\log(eN/m)}{m}\Big]^{1/p}\quad,\quad 1\leq m\leq N\,,
\end{equation}
with $``s_m \in \{d_m,e_m\}"$. The corresponding matching bounds for entropy numbers in the more general situation $\id:\ell_p^N \to \ell_q^N$ where $0<p< q \leq \infty$ is nowadays well-known but has a long history. Let us mainly refer to Sch\"utt \cite{Schu84}, Edmunds, Triebel \cite{EdTr96} and K\"uhn \cite{Ku01}. For the complete history of this result we refer to \cite[Rem.\ 3]{MaUl20}. The corresponding result for Kolmogorov numbers in case $2\leq p< \infty$ can be found in \cite[Thm.\ 1.1]{FoPaRaUl10}, where the sharp dual version in terms of Gelfand numbers is proved. This result in combination with the proof in \cite[Thm.\ 3.19]{Vy06} gives for all $2\leq p \leq \infty$, $r>1/p$
\begin{equation}\label{d-1}
	s_m(\operatorname{I}:\bB^r_{p,\infty} \to L_\infty)
	\leq s_m(\operatorname{I}:\bB^r_{p,\infty} \to \bB^0_{\infty,1}) \lesssim m^{-r} (\log m)^{(d-1)(r+1)}\,.
\end{equation}
It turns out that in case $d=2$ we recover the small smoothness result in Theorem \ref{theorem:main2b} as well as Belinskii's ``large smoothness'' result in \cite{TrBe04}, 11.3.5. In case $d=2$ the above result is sharp, see the discussion in Remark \ref{rem:d=2}. In addition, the result gives an indication that the $\log\log m$ term in Theorem \ref{theorem:main3} is probably not needed. In case of small smoothness and $d>2$ the result is worse than our result in Theorem \ref{theorem:main2b} and Theorem \ref{theorem:main3}.

Additionally, in \cite{Vy06} the author pointed out some gaps between upper and lower bounds in a certain range of small smoothness. This was the starting point of the recent paper Mayer, Ullrich \cite[Cor.\ 23, (iii)]{MaUl20}, where the sharp behavior \begin{equation}\label{res:MaUl}
   e_m(\operatorname{I}:\bB^r_{p,q} \to L_\infty) \asymp m^{-r}
\end{equation}
is shown in case $2 < p \leq \infty$, $0<q\leq 2/3$, and $1/p<r \leq 1/2$ for all dimensions $d$. The proof relies on a refinement of \eqref{f39} for mixed $\ell_q^n(\ell_p^N)$-norms, see \cite[Thm.\ 13]{MaUl20}. Roughly speaking, refining the spaces $\bH^r_p$ by decreasing the third parameter $q$, see Remark \ref{Fourier2}, allows us to get rid of the logarithmic term. A combination of the technique in \cite{MaUl20} with the technique used in this paper may allow to extend the range of parameters for the result \eqref{res:MaUl}. A corresponding result for Kolmogorov numbers is not known. However, a similar phenomenon occurs for the $\sigma_m$ numbers associated to a wavelet type dictionary (see below). Note that the space  $\bB^r_{p,q}$ is a quasi-Banach space.

\paragraph{Wavelet type dictionaries.} In the context of function spaces with mixed smoothness not only best $m$-term trigonometric approximation has been considered. Also hyperbolic wavelet type dictionaries $\mathcal{D}$, as defined in \eqref{f50},
gained substantial interest, see for instance \cite{TeGr00} or \cite[Sect.\ 7.2]{DuTeUl19} and the references therein. It turned out that the order of decay of the corresponding error quantities (modify the definition of $(Y_n)_n$ in \eqref{f31} accordingly) is often substantially better than for the  trigonometric system, see \cite[Sect.\ 7]{DuTeUl19} and the references therein. In fact, the gain is not only in the logarithmic term but sometimes also in the main rate. This is certainly not the case in our setting. A reasonable question would be: Does a wavelet system perform comparably well with respect to the decay of the associated $\sigma_m$ in the $L_\infty$ norm? Note that a fundamental difference between wavelets and the trigonometric system is the lack of a universal $L_\infty$-bound of the $L_2$-normalized wavelet system. As a consequence, some of the techniques used by Belinskii, see for instance \cite[11.2.5]{TrBe04}, can not be directly adapted to wavelets. A strong indication that wavelet dictionaries may not perform worse than the trigonometric system is the following observation. In \cite[Thm.\ 6.15]{By18} it is proved that for $2<p\leq \infty$, $0<q\leq 2/3$, and $1/p<r\leq 1/2$ we have
$$
   \sigma_m(\operatorname{I}:\bB^r_{p,q}(\T^d) \to \Linfty; \mathcal{D}) \lesssim m^{-r}\,.
$$
This result is sharp for all dimensions $d$ if we use the tensorized Faber Schauder system as dictionary $\mathcal{D}$.

%

\paragraph{Sampling recovery.} We introduce the notion of sampling numbers of an operator $T:\bF \to \bG$ between two Banach spaces $\bF$ and $\bG$ of functions on $D$. We assume that point evaluations are linear functionals on $\bF$. This would be the case if $\bF$ is continuously embedded into $\mathcal{C}(D)$, the space of continuous functions on $D$. Let us define the $m$-th sampling numbers of an operator $T \in \mathcal{L}(\bF,\bG)$ as follows
$$
	\varrho_m(T:\bF \to \bG) := \inf\limits_{\bx^1,...,\bx^n \in D} \; \inf\limits_{\substack{\varphi:\C^n \to \bG\\\text{linear}}} \; \sup\limits_{\|f\|_{\bF}\leq 1}\|Tf-\varphi(f(\bx^1),...,f(\bx^n))\|_{\bG}\,.
$$

In many cases the embedding $\bF \hookrightarrow \bG$ and the corresponding embedding operator $\operatorname{I}:\bF \to \bG$ is considered. A particular situation is the case when $\bG = L_2(D)$. In this situation it has been proven in \cite{temlyakov2020optimal} that there are two positive absolute constants $b,B>0$ such that
\begin{equation}\label{rel}
	\varrho_{bm}(\operatorname{I}:\bF \to L_2(D)) \leq Bd_m(\operatorname{I}:\bF \to L_\infty(D))\quad,\quad m\in \N\,.
\end{equation}
In case that $\bF$ represents a reproducing kernel Hilbert space $H(K)$ embedded into $L_2(D)$ we even know that, see \cite{NaSchUl20} and \cite{KrUl19}
$$
\varrho_m(\operatorname{I}:H(K) \to L_2)^2 \leq \frac{c_1\log(m)}{m}\sum\limits_{k=\lfloor c_2m \rfloor}^{\infty} d_k(\operatorname{I}:H(K) \to L_2(D))^2\quad,\quad m\in \N\,,
$$
with (precisely given) absolute constants $c_1,c_2>0$. Similar results have been recently  established for non-Hilbert function spaces, see \cite{KrUl20}. However, in both these settings the square summability of the corresponding Kolmogorov numbers
$d_k(\operatorname{I}:\bF \to L_2)$ is crucial. When considering the target space $\bG = L_\infty(D)$ we refer to the recent results \cite{PoUl21}.

Let us also mention a recent lower bound in the Hilbert space situation. The authors of \cite{HiKrNoVy21} showed that there exists a RKHS $H(K)$ with square summable Kolmogorov numbers $d_k(\operatorname{I}:H(K) \to L_2) \asymp k^{-1/2}(\log k)^{-\beta}$, $\beta>1/2$, such that $\varrho_m(I) \gtrsim m^{-1/2}(\log m)^{-\beta +1/2}$. This shows that already in the full Hilbert spaces setting a logarithmic gap between sampling and Kolmogorov numbers is possible. The slow decay of the Kolmogorov numbers could be interpreted as a certain ``small smoothness''. Note that in case of only the target space is $L_2$ (and the source space is non-Hilbert) such a gap has been observed earlier, see \eqref{f29} and \eqref{f29b} below.

In our ``small smoothness setting'' the decay rate of the Kolmogorov numbers (in $L_2$) is allowed to be strictly smaller than $1/2$. We are not restricted to the square summability of those. Hence, results on Kolmogorov numbers in the uniform norm together with \eqref{rel} serve as a powerful tool to investigate the  sampling recovery problem in $L_2$ for the  case of small smoothness. From \eqref{rel} together with Theorem \ref{theorem:main2b} we obtain in case $1/p<r<1/2$
\begin{equation}\label{f25}
	\varrho_m(\operatorname{I}:\Hpr \to \Ltwo) \lesssim m^{-r}(\log m)^{d-1+r}\,.
\end{equation}
In addition, the endpoint result from Theorem \ref{theorem:main3} have direct counterparts for sampling numbers. Let us point out that the so far best-known upper bounds in the above situation have been obtained by the use of sparse grid (Smolyak) recovery algorithms, see \cite{Tem93}, \cite{SiUl07}, \cite{Du11, Du16}, resulting in
\begin{equation} \label{f26}
   	\varrho_m(\operatorname{I}:\Hpr \to \Ltwo) \lesssim m^{-r}(\log m)^{(d-1)(1+r)}
\end{equation}
for all $r>1/p$. It is obvious, that \eqref{f25} improves on \eqref{f26} in case $d>2$ if $1/p<r < 1/2$. Note also that \cite[Thm.\ 5.1]{DuUl15} shows, when restricting to sparse grid methods, the bound in \eqref{f26} can not be improved. Hence, in the case of small smoothness sparse grid methods can not be optimal in the above situation when $d>2$.

As for $\operatorname{I}:\Wpr \to \Ltwo$ we obtain by \eqref{rel} and Theorem \ref{theorem:main1} the bound
\begin{equation}\label{f27}
	\varrho_m(\operatorname{I}:\Wpr \to \Ltwo) \lesssim m^{-r}(\log m)^{(d-1)(1-r)+r}
\end{equation}
if $1/p<r<1/2$. Clearly, the bound in Theorem \ref{theorem:main3} on the endpoint case also carries over to the sampling numbers. By the results in \cite{SiUl07} together with complex interpolation we find the bound
\begin{equation}\label{f28}
	\varrho_m(\operatorname{I}:\Wpr \to \Ltwo) \lesssim m^{-r}(\log m)^{(d-1)(1+\varepsilon)}\,,
\end{equation}
for any $\varepsilon>0$ in the case of small smoothness by using a sparse grid method. Clearly, \eqref{f27} improves on \eqref{f28}. However, here it is not clear whether the analysis for the sparse grid method can be improved or not. A valid lower bound comes from the embedding $\bB^r_{p,2} \hookrightarrow \bW^r_p$ since $p>2$. Hence, \cite[Thm.\ 5.1]{DuUl15} shows that any sparse grid method is asymptotically worse than $m^{-r}(\log m)^{(d-1)(r+1/2)}$. This yields in case
$1/4<r<1/2$ and large enough $d$ that the sparse grid methods can not be optimal since the bound in \eqref{f27} is better. However, the sampling method behind the bounds in this paper is highly non-constructive, whereas the sparse grid methods are constructive and can be implemented. From this point of view, our results show that there might exist further constructive methods which improve on the sparse grid methods regarding the asymptotic error decay.

Let us finally mention that the correct order of the quantities $\varrho_m$ in the above situation is still unknown. We improved on the upper bounds. A trivial source for lower bounds for $\varrho_m$ are the error quantities with respect to numerical integration. This in connection with the lower bounds in \cite{VTbookMA} and \cite{DuUl15} shows in case $1/p<r<1/2$ and $d\geq 2$
\begin{equation}\label{f29}
	\varrho_m(\operatorname{I}:\Wpr \to \Ltwo) \gtrsim m^{-r}(\log m)^{(d-1)/2} \gnapprox
	d_m(\operatorname{I}:\Wpr \to \Ltwo)
\end{equation}
and
\begin{equation}\label{f29b}
	\varrho_m(\operatorname{I}:\Hpr \to \Ltwo) \gtrsim m^{-r}(\log m)^{d-1} \gnapprox
	d_m(\operatorname{I}:\Hpr \to \Ltwo)\,.
\end{equation}
In fact, these bounds show that in the small smoothness range and $d\geq 2$ the problem of numerical integration (and hence the sampling recovery problem) is more ``difficult'' than the linear approximation problem (even in $L_p(\T^d)$). The corresponding worst-case errors decay slower by a $d$-dependent logarithm compared to the  corresponding Kolmogorov numbers, see \cite[Thms.\ 4.3.1, 4.3.10]{DuTeUl19}. This gap can get large if $d$ grows. An effect which has been observed earlier in \cite{HiNoVy09} for slowly decaying singular numbers in the Hilbert space setting. This can not happen if the corresponding Kolmogorov numbers (in $L_2$) decay fast enough (\cite{NaSchUl20, KrUl20}. Note that this is a multivariate effect since in dimension $d=1$ all the involved asymptotic characteristics decay as $m^{-r}$ as \eqref{d-1} shows (in case of $\sigma_m$ see \cite{Belinskii1998DecompositionTA}).

Furthermore, comparing \eqref{f29b} to \eqref{f25}, we determined the sampling numbers in $L_2$ up to a rather small $d$-independent logarithmic gap with exponent $r<1/2$. We would like to emphasize that the case of mixed H\"older continuous functions is included here such that we obtain in case $0<r<1/2$
$$
  m^{-r}(\log m)^{d-1} \lesssim \varrho_m(\operatorname{I}:\bH^r_{\infty}(\T^d) \to \Ltwo) \lesssim m^{-r}(\log m)^{d-1+r}\,.
$$

On the other hand, comparing \eqref{f29} to \eqref{f27} the difference in the $\log$-exponent is again growing in $d$ in case $r<1/2$. If $r=1/2$ we are close to the lower bound coming from numerical integration.

Finally, as a consequence of \eqref{rel} (see also \cite{No86}), we obtain the same chain of inequalities for $d_m$ instead of $\varrho_m$ in \eqref{f29} and \eqref{f29b}. In other words, in case $d\geq 2$ and $1/p<r<1/2$ it holds
\begin{equation}\label{dm_final}
	d_m(\operatorname{I}:\Hpr \to \Ltwo)) = o(d_m(\operatorname{I}:\Hpr \to L_\infty(\T^d))).
\end{equation}
Note, that \eqref{dm_final} is not true in case $d=1$ since all involved quantities decay as $m^{-r}$. In case $d=2$ the relation \eqref{dm_final} is true if $r \neq 1/2$. In all other cases it is not known.

\nocite{Seeger2009,No86}


\paragraph{Acknowledgment.}
The first author was supported by the Russian Federation Government Grant N{\textsuperscript{\underline{o}}}14.W03.31.0031. T.\,U. would like to acknowledge support by the DFG Ul-403/2-1. T.\,U. also thanks Winfried Sickel and Thomas K\"uhn for several discussions on the topic. Last but not least T.\,U. would like to thank Kevin Matthes who implemented parts of the manuscript in \TeX.





\begin{thebibliography}{10}

	\bibitem{Belinskii1998DecompositionTA}
	E.~Belinskii.
	\newblock Decomposition theorems and approximation by a ``floating'' system of
	  exponentials.
	\newblock {\em Transactions of the American Mathematical Society}, 350:43--53,
	  1998.

	\bibitem{Bel98}
	E.~S. Belinsky.
	\newblock Estimates of entropy numbers and {G}aussian measures for classes of
	  functions with bounded mixed derivative.
	\newblock {\em J. Approx. Theory}, 93(1):114--127, 1998.

	\bibitem{BeLo76}
	J.~Bergh and J.~L\"{o}fstr\"{o}m.
	\newblock {\em Interpolation spaces. {A}n introduction}.
	\newblock Springer-Verlag, Berlin-New York, 1976.
	\newblock Grundlehren der Mathematischen Wissenschaften, No. 223.

	\bibitem{By18}
	G.~Byrenheid.
	\newblock {\em {Sparse representation of multivariate functions based on
	  discrete point evaluations}}.
	\newblock Dissertation, Institut f{\"u}r Numerische Simulation,
	  Universit{\"a}t Bonn, 2018.

	\bibitem{Ca81}
	B.~Carl.
	\newblock Entropy numbers, $s$-numbers, and eigenvalue problems.
	\newblock {\em J. Funct. Analysis}, 41:290--306, 1981.

	\bibitem{Du11}
	D.~D\~{u}ng.
	\newblock B-spline quasi-interpolant representations and sampling recovery of
	  functions with mixed smoothness.
	\newblock {\em J. Complexity}, 27(6):541--567, 2011.

	\bibitem{DuTeUl19}
	D.~{D}\~ung, V.~N. {T}emlyakov, and T.~{U}llrich.
	\newblock {\em {H}yperbolic {C}ross {A}pproximation}.
	\newblock Advanced Courses in Mathematics. CRM Barcelona.
	  Birkh\"auser/Springer, 2019.

	\bibitem{DuUl15}
	D.~D\~ung and T.~Ullrich.
	\newblock Lower bounds for the integration error for multivariate functions
	  with mixed smoothness and optimal {F}ibonacci cubature for functions on the
	  square.
	\newblock {\em Math. Nachr.}, 288(7):743--762, 2015.

	\bibitem{Du16}
	D.~D{\~u}ng.
	\newblock Sampling and cubature on sparse grids based on a {B}-spline quasi-
	  interpolation.
	\newblock {\em Found. Comput. Math.}, 16(5):1193--1240, 2016.

	\bibitem{DuLiKuLi98}
	T.~Dunker, T.~K\"{u}hn, M.~Lifshits, and W.~Linde.
	\newblock Metric entropy of the integration operator and small ball
	  probabilities for the {B}rownian sheet.
	\newblock {\em C. R. Acad. Sci. Paris S\'{e}r. I Math.}, 326(3):347--352, 1998.

	\bibitem{EdTr96}
	D.~E. Edmunds and H.~Triebel.
	\newblock {\em Function spaces, entropy numbers, differential operators},
	  volume 120 of {\em Cambridge Tracts in Mathematics}.
	\newblock Cambridge University Press, Cambridge, 1996.

	\bibitem{FoPaRaUl10}
	S.~Foucart, A.~Pajor, H.~Rauhut, and T.~Ullrich.
	\newblock The {G}elfand widths of {$\ell_p$}-balls for {$0<p\leq 1$}.
	\newblock {\em J. Complexity}, 26(6):629--640, 2010.

	\bibitem{HiKrNoVy21}
	A.~Hinrichs, D.~Krieg, E.~Novak, and J.~Vybiral.
	\newblock Lower bounds for integration and recovery in $L_2$.
	\newblock 2021.

	\bibitem{HiMaOeUl16}
	A.~Hinrichs, L.~Markhasin, J.~Oettershagen, and T.~Ullrich.
	\newblock Optimal quasi-{M}onte {C}arlo rules on order 2 digital nets for the
	  numerical integration of multivariate periodic functions.
	\newblock {\em Numer. Math.}, 134(1):163--196, 2016.

	\bibitem{HiNoVy09}
	A.~Hinrichs, E.~Novak, and J.~Vyb\'{i}ral.
	\newblock Linear information versus function evaluations for
	  {$L_2$}-approximation.
	\newblock {\em J. Approx. Theory}, 153(1):97--107, 2008.

	\bibitem{KrUl19}
	D.~Krieg and M.~Ullrich.
	\newblock Function Values Are Enough for {$L_2$}-Approximation.
	\newblock {\em Found. Comput. Math.}, 21(4):1141--1151, 2021.

	\bibitem{KrUl20}
	D.~Krieg and M.~Ullrich.
	\newblock Function values are enough for {$L_2$}-approximation: {P}art {II}.
	\newblock {\em J. Complexity}, 66:Paper No. 101569, 14, 2021.

	\bibitem{Ku01}
	T.~K\"{u}hn.
	\newblock A lower estimate for entropy numbers.
	\newblock {\em J. Approx. Theory}, 110(1):120--124, 2001.

	\bibitem{MaUl20}
	S.~Mayer and T.~Ullrich.
	\newblock Entropy numbers of finite dimensional mixed-norm balls and function
	  space embeddings with small mixed smoothness.
	\newblock {\em Constr. Approx.}, 53(2):249--279, 2021.

	\bibitem{NaSchUl20}
	N.~Nagel, M.~Sch\"afer, and T.~Ullrich.
	\newblock A new upper bound for sampling numbers.
	\newblock {\em Found. Comput. Math.}, {\tt https://doi.org/10.1007/s10208-021-09504-0}.

	\bibitem{No86}
	E.~Novak.
	\newblock Quadrature and widths.
	\newblock {\em J. Approx. Theory}, 47:195--202, 1986.

	\bibitem{PaTo86}
	A.~Pajor and N.~Tomczak-Jaegermann.
	\newblock Subspaces of small codimension of finite-dimensional {B}anach spaces.
	\newblock {\em Proc. Amer. Math. Soc.}, 97(4):637--642, 1986.

	\bibitem{Pietsch78}
	A.~Pietsch.
	\newblock {\em Operator ideals}.
	\newblock North-Holland, 1980.

	\bibitem{Pie81}
	A.~Pietsch.
	\newblock Approximation spaces.
	\newblock {\em J. Approx. Theory}, 32(2):115--134, 1981.

	\bibitem{Pie87}
	A.~Pietsch.
	\newblock {\em Eigenvalues and {$s$}-numbers}, volume~13 of {\em Cambridge
	  Studies in Advanced Mathematics}.
	\newblock Cambridge University Press, Cambridge, 1987.

	\bibitem{PoUl21}
	K.~Pozharska and T.~Ullrich.
	\newblock A note on sampling recovery of multivariate functions in the uniform
	  norm.
	\newblock {\em arXiv:2103.11124}, 2021.

	\bibitem{Rom19}
	A.~S. Romanyuk.
	\newblock Entropy numbers and widths for the {N}ikol'skij-{B}esov classes of
	  functions of many variables in the space $L_{\infty}$.
	\newblock {\em Analysis Math.}, 45(1):133--151, 2019.

	\bibitem{Schu84}
	C.~Sch\"{u}tt.
	\newblock Entropy numbers of diagonal operators between symmetric {B}anach
	  spaces.
	\newblock {\em J. Approx. Theory}, 40(2):121--128, 1984.

	\bibitem{Seeger2009}
	A.~Seeger and W.~Trebels.
	\newblock Low regularity classes and entropy numbers.
	\newblock {\em Archiv der Mathematik}, 92:147--157.

	\bibitem{SiUl07}
	W.~Sickel and T.~Ullrich.
	\newblock The {S}molyak algorithm, sampling on sparse grids and function spaces
	  of dominating mixed smoothness.
	\newblock {\em East J. Approx.}, 13(4):387--425, 2007.

	\bibitem{Tem93}
	V.~N. Temlyakov.
	\newblock {\em Approximation of periodic functions}.
	\newblock Computational Mathematics and Analysis Series. Nova Science
	  Publishers Inc., Commack, NY, 1993.

	\bibitem{TeGr00}
	V.~N. Temlyakov.
	\newblock Greedy algorithms with regard to multivariate systems with special
	  structure.
	\newblock {\em Constr. Approx.}, 16:399--425, 2000.

	\bibitem{Te15}
	V.~N. Temlyakov.
	\newblock Constructive sparse trigonometric approximation and other problems
	  for functions with mixed smoothness.
	\newblock {\em Matem. Sb.}, 206:131--160, 2015.

	\bibitem{VTbookMA}
	V.~N. Temlyakov.
	\newblock {\em Multivariate approximation}, volume~32 of {\em Cambridge
	  Monographs on Applied and Computational Mathematics}.
	\newblock Cambridge University Press, Cambridge, 2018.

	\bibitem{temlyakov2020optimal}
	V.~N. Temlyakov.
	\newblock On optimal recovery in {$L_2$}.
	\newblock {\em J. Complexity}, 65:Paper No. 101545, 11, 2021.

	\bibitem{Te20_5}
	V.~N. Temlyakov.
	\newblock Sampling discretization of integral norms of the hyperbolic cross
	  polynomials.
	\newblock {\em Tr. Mat. Inst. Steklova}, 312:282--293, 2021.

	\bibitem{TU20}
	V.~N. Temlyakov and T.~Ullrich.
	\newblock Bounds on {K}olmogorov widths and sampling recovery for classes with
	  small mixed smoothness.
	\newblock {\em J. Complexity}, 67:Paper No. 101575, 2021.

	\bibitem{TrBe04}
	R.~M. Trigub and E.~S. Bellinsky.
	\newblock {\em Fourier analysis and approximation of functions}.
	\newblock Kluwer Academic Publishers, Dordrecht, 2004.
	\newblock [Belinsky on front and back cover].

	\bibitem{UlUl16}
	M.~Ullrich and T.~Ullrich.
	\newblock The role of {F}rolov's cubature formula for functions with bounded
	  mixed derivative.
	\newblock {\em SIAM J. Numer. Anal.}, 54(2):969--993, 2016.

	\bibitem{Vy06}
	J.~Vyb{\'i}ral.
	\newblock Function spaces with dominating mixed smoothness.
	\newblock {\em Dissertationes Math. (Rozprawy Mat.)}, 436:73, 2006.

	\end{thebibliography}


\end{document}